\documentclass{ifacconf}

\usepackage{graphicx}     
\usepackage{natbib}
\usepackage{color}
\usepackage{courier}
\usepackage{algorithm}
\usepackage[noend]{algpseudocode}
\usepackage{amsmath}
\usepackage{amssymb}
\usepackage{amsfonts}
\usepackage{array}    

\usepackage{float}
\usepackage{tikz}
\usetikzlibrary{arrows,%
                plotmarks, calc}
\usetikzlibrary[positioning, calc, intersections, backgrounds]      

\usepackage{xspace}
\usepackage{soul}
\usepackage{mathtools}
\usepackage{makecell}
\usepackage{cite}
\usepackage{nicefrac}
\usepackage{float}
\usepackage{wrapfig}
\usepackage[utf8]{inputenc}

\newcommand{\ra}{\rightarrow}

\newtheorem{dfn}{Definition}

\newtheorem{exm}{Example}

\newcommand{\ie}{\unskip, i.\,e.,\xspace}
\newcommand{\eg}{\unskip, e.\,g.,\xspace}

\newcommand{\sut}{\text{s.\,t.\,}}
\newcommand{\wrt}{w.\,r.\,t.\xspace}

\newcommand{\N}{\ensuremath{\mathbb{N}}}

\newcommand{\R}{\ensuremath{\mathbb{R}}}

\newcommand{\X}{\ensuremath{\mathbb{X}}}
\newcommand{\Y}{\ensuremath{\mathbb{Y}}}
\newcommand{\F}{\ensuremath{\mathbb{F}}}
\newcommand{\U}{\ensuremath{\mathbb{U}}}

\newcommand{\sm}{\ensuremath{\setminus}}
\newcommand{\set}[1]{\ensuremath{\mathbb{#1}}}

\let\emptyset\varnothing

\newcommand{\cco}{\ensuremath{\overline{\text{co}}}}

\newcommand{\eps}{\ensuremath{\varepsilon}}

\newcommand{\pdiff}[2]{ \ensuremath{ \frac{\partial {#1}}{\partial {#2}} } }

\newcommand{\ball}{\ensuremath{\mathcal B}}
\newcommand{\D}{\ensuremath{\mathcal{D}}}		
\DeclareMathOperator*{\argmin}{arg\,min}

\DeclareMathOperator*{\arginf}{arg\,inf}

\newcommand{\nrm}[1]{\left\lVert#1\right\rVert}

\newcommand{\abs}[1]{\left\lvert#1\right\rvert}
\newcommand{\scal}[1]{\left\langle#1\right\rangle}

\definecolor{dgreen}{rgb}{0.0, 0.5, 0.0}

\newcommand{\spc}{\ensuremath{\,\,}}	

\makeatletter
\newcommand{\subalign}[1]{%
	\vcenter{%
		\Let@ \restore@math@cr \default@tag
		\baselineskip\fontdimen10 \scriptfont\tw@
		\advance\baselineskip\fontdimen12 \scriptfont\tw@
		\lineskip\thr@@\fontdimen8 \scriptfont\thr@@
		\lineskiplimit\lineskip
		\ialign{\hfil$\m@th\scriptstyle##$&$\m@th\scriptstyle{}##$\crcr
			#1\crcr
		}%
	}
}

\newcolumntype{L}[1]{>{\raggedright\let\newline\\\arraybackslash\hspace{0pt}}m{#1}}
\newcolumntype{C}[1]{>{\centering\let\newline\\\arraybackslash\hspace{0pt}}m{#1}}
\newcolumntype{R}[1]{>{\raggedleft\let\newline\\\arraybackslash\hspace{0pt}}m{#1}}

\usepackage[inline]{enumitem}
       
\begin{document}

\setlength{\abovedisplayskip}{0pt}
\setlength{\belowdisplayskip}{0pt}
\setlength{\abovedisplayshortskip}{0pt}
\setlength{\belowdisplayshortskip}{0pt}	
	
\begin{frontmatter}

\title{Nonsmooth stabilization and its computational aspects}


\author[TUC]{Pavel Osinenko, Patrick Schmidt, Stefan Streif}

\address[TUC]{Technische Universit\"at Chemnitz, Automatic Control and System Dynamics Laboratory, Germany (e-mail: p.osinenko@gmail.com, \{patrick.schmidt, stefan.streif\}@etit.tu-chemnitz.de)}

\begin{abstract}
This work has the goal of briefly surveying some key stabilization techniques for general nonlinear systems, for which, as it is well known, a smooth control Lyapunov function may fail to exist.
A general overview of the situation with smooth and nonsmooth stabilization is provided, followed by a concise summary of basic tools and techniques, including general stabilization, sliding-mode control and nonsmooth backstepping.
Their presentation is accompanied with examples.
The survey is concluded with some remarks on computational aspects related to determination of sampling times and control actions.
\end{abstract}

\begin{keyword}
Stabilization, nonsmooth analysis, control Lyapunov function
\end{keyword}

\end{frontmatter}




\section{Introduction} \label{sec:intro}

Nonsmooth tools play central role in nonlinear system stabilization theory as the smooth ones are prone to limitations posed by the celebrated work of \citet{Brockett1983-stabilization}.
He demonstrated that even such simple systems, as a three-wheel robot, do not admit a smooth control Lyapunov function (CLF).
A particular consequence of this fact is that there can be no continuous control law, which depends only on the system's state and which parks the robot into the desired position.
The alternatives are either to consider a \emph{time-varying} (or dynamical) control law or to give up the continuity condition.
The former approach received great attention in the 80s and 90s \citep{Aeyels1985-stabilization,Kawski1989-stabilization-plane,Coron1991-smooth-stabilization,Samson1991-stabilization,Pomet1992-dyn-stabilization,Coron1992-stabilization,Coron1993-stab-time-var,Coron1994-stabilization,Coron1995-stabilization,Khaneja1999-dynamic-feedback,Morin1999-design}.
Although the design of time-varying control laws may happen to be somewhat involved, the great advantage of this approach is that it requires the usual analysis tools and the closed-loop system trajectory is a classical Carathéodory solution which enjoys uniqueness properties.
Contrary to this approach, if the continuity condition of the control law is omitted, care must be taken when defining what a system trajectory actually is.
That is where alternative solutions to the respective initial value problems with discontinuous right-hand side come into play.
One of the most well-known is the Filippov solution \citep{Filippov2013-discont-dyn-sys}.
In brief, it is an absolutely continuous function that satisfies the said initial value problem with the right-hand side interpreted in the sense of a \emph{differential inclusion} of the kind $\dot x \in F(x,t)$, where $F$ is a \emph{set-valued map}.
Sliding-mode control (SMC) makes wide use of Filippov solutions \citep{Slotine1991-nonlin-ctrl,Young1996-slid-mode,Fridman1996-slid-mode,Fridland1999-slid-mode,Perruquetti2002-slid-mode}.
A particular drawback of this kind of a solution in the context of stabilization is that the very same Brockett's conditions apply just as if one were to limit to the classical solutions \citep{Coron1994-stabilization,Ryan1994-stabilization}.
Further solutions include the ones in the sense of Hermes, Krasovskii, Sentis etc. whose good overview was done by \citet{Cortes2008-discont-dyn-sys}.

Of particular interest for this survey is the setting of \emph{sample-and-hold} (S\&H) solutions which are very simple to interpret.
For the original nonlinear system
\begin{equation} \label{eqn:sys}
	\tag{Sys}
	\dot x  = f(x, u), x \in \R^n, u \in \R^m	
\end{equation}
if one designs a discontinuous static control law $\kappa: \R^n \ra \R^m$, instead of applying $\kappa$ literally, one can consider its \emph{digital implementation} with a sampling time $\delta > 0$ in the form
\begin{equation} \label{eqn:sys-SH}
	\tag{Sys-SH} 
	\begin{aligned}
		& \dot x = f(x, u_k), \\
		& t \in [ k \delta, (k + 1) \delta], u_k \equiv \kappa(x (k \delta)), k \in \N.
	\end{aligned}
\end{equation}
In this case, if $f$ is locally Lipschitz \wrt $x$, classical trajectories of \eqref{eqn:sys-SH} exist (at least locally) and are unique, since $f$ is also measurable in $t$, precisely due to the sample nature of the control.
The only property one has to give up in general when stabilizing \eqref{eqn:sys} in the S\&H mode, is asymptotic stabilizability of the closed loop.
Instead, one has \emph{practical stabilizability} in the sense of
\begin{dfn}[Practical stabilizability]
	\label{dfn:pract-stab-approx-opt}
	A control law $\kappa$ is said to practically stabilize \eqref{eqn:sys} in the sense of S\&H \ie \eqref{eqn:sys-SH} if, for any data $R > r > 0$, there is a sufficiently small $\delta > 0$ such that any closed-loop trajectory $x(t), t \ge 0, x(0) \in \ball_R$ is bounded and enters and stays in the ball $\ball_r$ within a time $T$ depending uniformly on $R, r$. 
\end{dfn} 
So, by properly selecting the sampling time, one can achieve any desirable precision of stabilization.
Thus, S\&H scenario of stabilization can be justified from both the implementation and usefulness sides.

\begin{table}[h]
	\caption{Overview of stabilization methods}
	\begin{tabular}{|L{0.2\columnwidth}|L{0.2\columnwidth}|L{0.2\columnwidth}|L{0.2\columnwidth}|}
		\hline		
		& \textbf{Cont., static} & \textbf{Cont., time-varying} & \textbf{Discont., static}
		\tabularnewline
		\hline	
		pros & Classical trajectories, rel. simple design & Classical trajectories & Rel. simple design
		\tabularnewline
		\hline	
		cons & Application limited & Sophisticated design & Nonstandard trajectories (however, S\&H)
		\tabularnewline 
		\hline
	\end{tabular}
	\label{tab:stabilization-methods}
\end{table}

A brief summary of the discussed stabilization methods can be found in Table \ref{tab:stabilization-methods}.
The focus of this survey is set to discontinuous static control laws in their S\&H realization due to the aforementioned usefulness and practicability.
In general, as said earlier, a nonlinear system does not admit a smooth CLF, so one has to widen the perspective to include nonsmooth tools, such as nonsmooth CLFs, generalized derivatives, subgradients etc., whose brief overview is given in the next section.
Particular stabilization techniques are surveyed in Section \ref{sec:stabilization}.
In Section \ref{sec:backstepping}, we overview a very powerful technique of nonsmooth backstepping with application to dynamically actuated three-wheel robot and dynamical Artstein's circles.
A brief survey of SMC in the S\&H setting is given in Section \ref{sec:SMC}.
The work is concluded with a discussion on computational aspects of stabilization techniques and a related case study in Section \ref{sec:comp-aspects}.

In the following, $\nrm x$ describes the Euclidean norm of $x$ and $\cco (\X)$ is defined as the closure of the convex hull of a set $\X$.
Furthermore, $\ball_R(x)$ denotes a ball with radius $R$ at $x$ \ie $\ball_R(x) := \{ x \in \R^n: \nrm x \leq R \}$ and $\ball_R$ means the same with $x = 0$.
Finally, $\rightrightarrows$ denotes a set-valued mapping.

\section{Basic nonsmooth tools} \label{sec:nonsmooth-basics}

Perhaps, the most central difference between the machinery of continuous and discontinuous stabilization lies at the level of a CLF.
As mentioned earlier, a general nonlinear system does not admit a smooth CLF, so one is forced to consider nonsmooth alternatives.
For this sake, let us first consider
\begin{dfn}[Lower directional generalized derivative]
	\label{dfn:dini-der}
	\textcolor{white}{For} For a locally Lipschitz function $V: \R^n \ra \R$, the function $\D_{\bullet} V(\bullet): \R^n \times \R^n \ra \R$ defined as
	\begin{align} \label{eqn:Dini-der}
		\begin{split}
			\D_{\vartheta} V(x) & \triangleq \liminf_{\mu \ra 0^+} \frac{ V(x + \mu \vartheta) - V(x)}{\mu}. 
			\end{split}
		\tag{Der}
	\end{align}
	is called \textit{lower directional generalized derivative} (further, just ``LDGD'').	
\end{dfn}

\begin{exm}
	\label{exm:Dini-der}
	Consider a function
	\[
		V(x) = \begin{cases}
			x^2 \sin \frac 1 x, & x < 0, \\
			0, & x = 0, \\
			2x^2 \sin \frac 1 x, & x > 0,
		\end{cases}
	\]
	LDGD at zero along one is $-2$.
	\begin{figure}[H]
		\centering \includegraphics{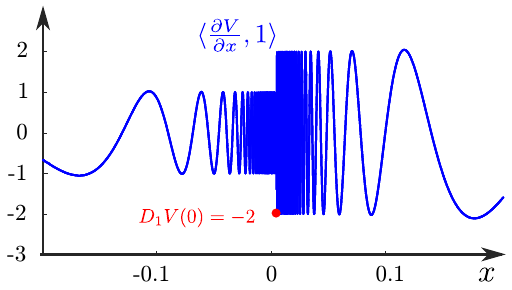}
		\caption{Graph of the derivative of $V$}
	\end{figure}
\end{exm}

Using the introduced LDGD, we can now consider one type of a nonsmooth CLF:
\begin{dfn}[CLF in LDGD sense]
	\label{dfn:dini-CLF}
	For the system \eqref{eqn:sys}, a locally Lipschitz, proper and positive-definite function $V: \R^n \ra \R$ is called \textit{CLF in LDGD sense} if there exists a continuous function $w: \R^n \ra \R, x \ne 0 \implies w(x) > 0$ satisfying a \emph{decay condition}: for any compact set $\X \subseteq \R^n$, there exists a compact set $\U_{\X} \subseteq \U$ such that
	\begin{equation}
		\label{eqn:Dini-decay-cond}
			\forall x \in \X \inf_{\vartheta \in \cco(f(x, \U_{\X}))} \; \D_{\vartheta} V(x) \le -w(x).
		\tag{Dec}
	\end{equation}
	The function $w$ is also called \textit{decay function}.	
\end{dfn}


The condition \eqref{eqn:Dini-decay-cond} effectively means that $V$ is an \textit{upper minimax solution} of the Hamilton-Jacobi equation
\begin{equation}
	\label{eqn:HJ}
	\tag{HJ}
	\inf_{u \in \U_{\X}}\scal{\nabla V, f(x,u)} + w(x) = 0
\end{equation}
on the respective domain \citep{Camilli2008-CLF-by-Zubov-mtd,Subbotin2013-gen-sol-PDE}.
There exist several techniques of practical stabilization of \eqref{eqn:sys} using an LDGD CLF in the sense of Definition \ref{dfn:dini-CLF}, some of which are reviewed in Section \ref{sec:stabilization}.
Here, it has to be clarified what ``$\nabla V$'' means in the nonsmooth setting.
There are several substitutes for gradients with different contexts, some of which are overviewed here. 
First of all, $\nabla V$ for a smooth function $V$ is a unique vector at each point.
If $V$ is nonsmooth, there is no unique vector which describes a descent direction of $V$, so one speaks of a set of those, summarizing them in a \textit{subdifferential}.
Here is the first such subdifferential, which is useful in practical stabilization \citep{Clarke1997-stabilization}:

\begin{dfn}[Proximal subdifferential]
	\label{dfn:prox-subdiff}
	A vector $\zeta \in \R^n$ is called \textit{proximal subgradient} of $V$ at $x \in \R^n$ if there exists a ball $\ball_r(x)$ and $\sigma > 0$ \sut
	\begin{equation}\label{eqn:prox-subgrad}
		\tag{Prox}
		\forall y \in \ball_r(x) \spc V(y) \geq V(x) + \scal{\zeta, y - x} - \sigma \nrm{y - x}^2. 
	\end{equation}
	The set of all such vectors is called \textit{proximal subdifferential} and is denoted by $\partial_P V(x)$.
%
\end{dfn}

It follows straight from the definition of an LDGD, that for any vector $\vartheta$ and any proximal subgradient $\zeta$, it holds that 
\begin{equation} \label{eqn:subgradient-zeta-theta}
	\scal{\zeta, \vartheta} \le \D_\vartheta V(x).
\end{equation}
Therefore, a decay condition in the spirit of \eqref{eqn:Dini-decay-cond} can be formulated as
\begin{equation}\label{eqn:prox-CLF}
	\tag{PDec}
	\forall \zeta \in \partial_P V(x) \spc \inf_{u \in \U_\X} \scal{\zeta, f(x,u)} \le - w(x).
\end{equation}
The condition \eqref{eqn:prox-CLF} means that $V$ is a \textit{proximal supersolution} of \eqref{eqn:HJ} \citep{Clarke1995-sys-traj}, or, equivalently, a \textit{viscosity supersolution} thereof \citep{Crandall1983-viscosity-HJB}.

Other subdifferentials exist, for instance:
\begin{dfn}[Limiting subdifferential]
	\label{dfn:limit-subdiff}
	Let $V: \R^n \ra \R$ be a function. 
	The set 
	\begin{equation} \label{eqn:limit-subdiff}
		\begin{aligned}
			\partial_L V(x) := \{ &\zeta \in \R^n: \zeta = \text{w-lim } \zeta_i, \\
			& \zeta_i \in \partial_P V(x_i), x_i \ra x, V(x_i) \ra V(x)\}
		\end{aligned}	
	\end{equation}
	is called \textit{limiting subdifferential} of $V$ at $x$, where w-lim is the weak limit.
\end{dfn}

Here, a sequence $\{ \zeta_i \}_{i = 1, 2, \ldots} \subset \R^n$ is said to converge weakly to $\zeta \in \R^n$, if $\scal{ \zeta_i, \theta } \ra \scal{ \zeta, \theta }$ for all $\theta \in \R^n$.
Such limiting constructions are used \eg in the nonsmooth practical stability analysis of SMC \citep{Clarke2009-slid-mode-stab}.
An overview is presented in Section \ref{sec:SMC}.

\begin{dfn}[Clarke subdifferential]
	\label{dfn:clarke-subdiff}
	Let $\partial_L V(x)$ be the limiting subdifferential of $V$ at $x$.
	Then, the \textit{Clarke subdifferential} $\partial_C V(x)$ is defined as the closed convex hull of \eqref{eqn:limit-subdiff}, i.e.
	\begin{equation}
	\partial_C V(x) = \cco \left( \partial_L V(x) \right).
	\end{equation}
\end{dfn}

\begin{exm}
	Consider $V(x) = \begin{cases}
	x^2, & \text{if } x < 1 \\
	(x-1)^2 + 1, & \text{else}	
	\end{cases}$. 
	The subdifferentials defined in Definitions \ref{dfn:prox-subdiff}, \ref{dfn:limit-subdiff} and \ref{dfn:clarke-subdiff} are given as
	\begin{align*}
	&\partial_P V(1) = \emptyset, \spc 
	\partial_L V(1) = \left\{ \begin{pmatrix} 1 \\ 0 \end{pmatrix}, \begin{pmatrix} 1 \\ 2 \end{pmatrix} \right\}, \\
	&\partial_C V(1) = \cco \left( \left\{ \begin{pmatrix} 1 \\ 0 \end{pmatrix}, \begin{pmatrix} 1 \\ 2 \end{pmatrix} \right\} \right).	
	\end{align*}
		\begin{figure}[H]
		\centering \includegraphics{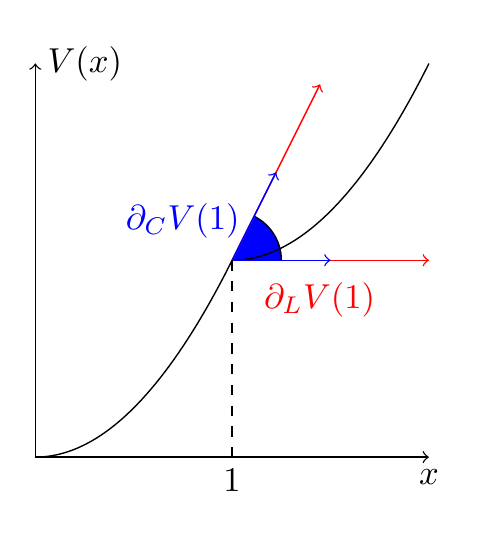}
		\caption{Graph of $V$ and the three subdifferentials}
	\end{figure}	
\end{exm}

An important property of many practical stabilizing techniques is \textit{semiconcavity}.
In fact, semiconcavity is a ubiquitous property of CLFs \citep{Clarke2011-discont-stabilization, Cannarsa2004-semiconcave}.
\begin{dfn}[Locally semiconcave function]
	\label{dfn:semiconc-fcn} 
	A function $V: \X \ra \R$ with $\X \subset \R^n$ being convex, is called \textit{locally semiconcave with linear modulus}, if it is continuous in $\X$ and there exists $C \geq 0$ such that the following inequality holds for all $x, y \in \X$
	\begin{equation}\label{eqn:semiconc}
		\tag{SemiConc}
		V(x) + V(y) - 2 V \left( \frac{x + y}{2} \right) \leq C \nrm{x - y}^2.
	\end{equation}
\end{dfn}

The following theorem states that any locally semiconcave function can be represented as the infimum of a family of $\mathcal C^2$ functions.
The proof can be found in \citep{Cannarsa2004-semiconcave}.

\begin{thm}
	\label{thm:semiconc-minimum-c2}
	Let $V: \X \ra \R, \X \subseteq \R^n$ be a locally semiconcave function with linear modulus according to Definition \ref{dfn:semiconc-fcn}.
	Then, for each compact subset $\set K \subset \X$, there exists a compact set $\Theta \subset \R^{2n}$ and a continuous function $F: \set K \times \Theta \ra \R$, \sut $F(\bullet; \theta)$ is $\mathcal C^2$ for any $\theta \in \Theta$, the gradients $\nabla_x F(\bullet; \theta)$ are equicontinuous, and $V$ can be expressed as
	\begin{equation} \label{eqn:marginal-fcn}
		V(x) = \min_{\theta \in \Theta} F(x; \theta),
	\end{equation}	
	for all $x \in \set K$.
	If there exists such a representation of $V$, then $V$ is called \textit{marginal function} \citep{Cannarsa2004-semiconcave}.
\end{thm}

Theorem \ref{thm:semiconc-minimum-c2} will be useful in presenting a particular technique of nonsmooth backstepping in Section \ref{sec:backstepping}.
Marginal functions are also used to define
\begin{dfn}[$F$-disassembled subdifferential]
	\label{dfn:disass-diff}
	\textcolor{white}{Let} Let $V: \X \ra \R, \spc \X \subseteq \R^n$ be a locally semiconcave function and let $F: \set K \times \Theta \ra \R$, where $\set K \subset \R^n$ and $\Theta \subseteq \R^{2n}$ are compact sets.
	The set-valued map $\partial_D^F V: \X \rightrightarrows \R^n$ defined as
	\begin{equation} \label{eqn:disass-diff}
		\partial_D^F V(x) \triangleq \left\{ \pdiff{F(x; \theta)}{x}: \theta \in \argmin_{\theta \in \Theta} F(x;\theta) \right\}
	\end{equation}
	is called \textit{$F$-disassembled subdifferential}.
	A single element of $\partial_D^F V(x)$ is called \textit{$F$-disassembled subgradient}.
\end{dfn}

Such a subdifferential was used in \citep{Cannarsa2004-semiconcave} and was named in \citep{Nakamura2013-asymptotic}.

%
%
%

Using disassembled subdifferentials, the following type of a nonsmooth CLF can be introduced \citep{Nakamura2013-asymptotic}
\begin{dfn}[$F$-disassembled CLF]
	\label{dfn:disassem-clf}
	A proper, positive-definite, locally semiconcave function $V: \R^n \ra \R$ is called \textit{$F$-disassembled CLF} for \eqref{eqn:sys}, if there exists a continuous, positive definite function $w: \R^n \ra \R$ satisfying the following decay condition: for any compact subset $\X \subseteq \R^n$, there exists a compact set $\U_{\X} \subseteq \U$, such that
	\begin{equation}\label{eqn:dissassem-CLF}
		\tag{DisDec}
		\begin{aligned}
			\forall x \in \X \spc \exists \zeta \in \partial_D^F V(x): \spc 
			&\min_{u \in \U_{\X}} \scal{\zeta, f(x, u)} \leq -w(x).
		\end{aligned}
	\end{equation}
\end{dfn}

An immediate relation between $F$-disassembled and proximal subdifferentials can be stated in the following lemma.
\begin{lem}
	\label{lem:prox-disass-subgr}
	Let $\X \subseteq \R^n$ be open and let $\Theta \subset \R^{2 n}$ be compact.  
	Let $V$ be a locally semiconcave function given as \eqref{eqn:marginal-fcn} for $F$ continuous in $\X \times \Theta$.
	If $\xi$ is a proximal subgradient of $V$ at $x$, then it is also an $F$-disassembled subgradient of $V$ at $x$ \ie 
	\begin{equation} \label{eqn:relation-disass-prox-subgr}
		\partial_P V(x) \subseteq \partial_L V(x) \subseteq \partial_D^F V(x).
	\end{equation}
\end{lem}

\begin{pf}
	The first inclusion follows from Definition \ref{dfn:limit-subdiff}.
	The second one is a part of the proof of Theorem 3.4.4 in \citep{Cannarsa2004-semiconcave}.
	\raggedleft \hfill $\blacksquare$ 
\end{pf}

The following lemma shows the relation between the two different kinds of CLFs.

\begin{lem}
	\label{lem:Prox-CLF-Disass-CLF}
	Consider \eqref{eqn:sys} and a related $F$-disassembled CLF $V: \R^n \ra \R$ with an arbitrary $F$.
	Then, $V$ is also a CLF in LDGD sense.
\end{lem}

\begin{pf}
	Choose a function $F: \X \times \Theta \ra \R$, $F \in \mathcal C^2(\X \times \Theta)$ satisfying \eqref{eqn:marginal-fcn} and let $\partial_D^F V(x)$ be the corresponding $F$-disassembled differential.
	Define for all $x \in \X$
	\begin{align*}
		\Phi (x, \vartheta, \eps; \theta) := \frac{F(x + \eps \vartheta; \theta) - F(x; \theta)}{\eps}
	\end{align*}
	Obviously, 
	\begin{equation} \label{eqn:Dini-lim-liminf}
		\begin{aligned}
			&\liminf_{\eps \ra 0^+} \Phi (x, \vartheta, \eps; \theta) \\
			&\le \lim_{\eps \ra 0} \Phi (x, \vartheta, \eps; \theta) = \scal{\pdiff{F(x; \theta)}{x}, \vartheta},
		\end{aligned}
	\end{equation}	
	holds for all $x \in \X$, which follows from the definition of $\liminf$.
	Since $\Phi$ is differentiable w.r.t. $x$ and $\theta$, the equality in \eqref{eqn:Dini-lim-liminf} holds for all $\theta \in \Theta$.
	In particular, for $\theta^\star \in \argmin_{\theta \in \Theta} F(x; \theta)$ it holds that
	\begin{equation}
		\begin{aligned}
			\lim_{\eps \ra 0^+} \Phi(x, \vartheta, \eps; \theta^\star) 
			= \scal{\left.\pdiff{F(x; \theta)}{x} \right \vert_{\theta = \theta^\star}, \vartheta} = \scal{\zeta, \vartheta}
		\end {aligned}
	\end{equation}
	for some $\zeta \in \partial_D^{F} V(x)$ corresponding to $\theta^\star$ at $x$.
	In turn, one obtains using \eqref{eqn:subgradient-zeta-theta}
	\begin{equation} \label{eqn:Disass-CLF-Dini-CLF}
		\begin{aligned}
			\lim_{\eps \ra 0^+} \Phi (x, \vartheta, \eps; \theta^\star) = \scal{\zeta, \vartheta} \le \D_\vartheta V(x).
		\end {aligned}
	\end{equation}
	Since $\lim_{\eps \ra 0^+} \Phi (x, \vartheta, \eps; \theta^\star)$ is independent of the choice of $F$, \eqref{eqn:dissassem-CLF} $\implies$ \eqref{eqn:Dini-decay-cond} holds with \eqref{eqn:Disass-CLF-Dini-CLF}.
	And since $F$ was chosen arbitrary, \eqref{eqn:Disass-CLF-Dini-CLF} holds for all $F$ that satisfy \eqref{eqn:marginal-fcn}.
	\raggedleft \hfill $\blacksquare$
\end{pf}

Before proceeding to concrete stabilizing techniques, the \textit{inf-convolution} (InfC) should be recalled \citep{Clarke2008-nonsmooth-analys}:
\begin{dfn}[Inf-convolution]
	\label{dfn:InfC}
	Let $V: \X \ra \R, \X \subseteq \R^n$. 
	For $\alpha \in (0,1)$, the \textit{inf-convolution} of $V$ at $x$ is defined by
	\begin{equation} \label{eqn:InfC}
		\tag{InfC}
		V_\alpha(x) \triangleq \inf_{y \in \R^n} \left\{ V(y) + \frac{1}{2 \alpha^2} \nrm{y - x}^2 \right\}. 	
	\end{equation}
\end{dfn}

In classical convex analysis, $V_\alpha$ is known as Yoreau-Mosida regularization of a (convex) function $V$. 
Furthermore, if $V$ is a lower semicontinuous function and bounded from below, then $V_\alpha$ is locally Lipschitz and an approximation of $V$ in the sense of $\lim_{\alpha \ra 0} V_\alpha(x) = V(x)$ \citep{Clarke1997-stabilization}.

\section{Stabilization techniques} \label{sec:stabilization}

In the following section, some stabilization techniques are presented.
Some of them are also described in \citep{Braun2017-SH-stabilization-Dini-aim}. 

\subsection{Steepest descent}

A steepest descent control law $\kappa: \R^n \ra \U \subseteq \R^m$ at $x \ne 0$ is computed via

\begin{equation} \label{eqn:steepest-descent}
	\kappa(x) \in \argmin_{u \in \U} \D_{f(x, u)} V(x).
\end{equation}

It is shown in \citep{Braun2017-SH-stabilization-Dini-aim} that $\kappa(x)$, computed by \eqref{eqn:steepest-descent} with a semiconcave LDGD CLF $V: \R^n \ra \R_{\geq 0}$, practically asymptotically stabilizes the origin of \eqref{eqn:sys}. 
The semiconcavity of $V$ is crucial to guarantee practical stabilizability by steepest descent \citep{Braun2017-SH-stabilization-Dini-aim}.

\subsection{Dini Aiming}

A control law at the state $x \in \R^n$ by Dini Aiming \citep{Kellett2000-Dini-aim} is computed in two steps based on a nondecreasing continuous function $\sigma: \R_{\geq 0} \ra \R_{\geq 0}$.

\begin{itemize}
	\item[1] Identify a direction $\vartheta^\star$ by means of minimizing the LDGD CLF $V$ over a neighborhood of $x$ for a given $r > 0$ \ie
		\begin{equation} \label{eqn:dini-aiming-direction}
			\vartheta^\star \in \argmin_{s \in \bar \ball_r (x)} V(s).
		\end{equation}
		
	\item[2] Compute an admissible control $\kappa(x) \in \U \cap \bar \ball_{\sigma(\nrm x + r)}$, $\U \cap \bar \ball_{\sigma(\nrm x + r)} = \{ u \in \U: \nrm u \leq \sigma(\nrm x + r) \}$ via
		\begin{equation} \label{eqn:dini-aiming-feeback}
			\kappa(x) = \argmin_{u \in \U \cap \bar \ball_{\sigma(\nrm x + r)} } \frac{ \scal{ x - \vartheta^\star, f(x, u) } }{ \nrm{ x - \vartheta^\star } }
		\end{equation}	
\end{itemize}

If the sampling time $\delta$ is chosen small enough (in accordance with $\sigma$ and $r$), the control computed in \eqref{eqn:dini-aiming-feeback} practically asymptotically stabilizes the origin of \eqref{eqn:sys} \citep{Kellett2004-Dini-aim, Braun2017-SH-stabilization-Dini-aim}.
Such a function $\sigma$ always exists, since $V$ is an LDGD CLF.
Furthermore, $\nrm x \ra 0 \implies \nrm u \ra 0$ has to hold.
In contrast to the steepest descent, $V$ does not necessarily need to be semiconcave.

\subsection{Optimization-based control}

In optimization-based control, the given LDGD is directly minimized over the set of admissible constant inputs \citep{Braun2017-SH-stabilization-Dini-aim} \ie
\begin{equation} \label{eqn:optim-based-feedb-optim-1}
	\min_{u \in \U} \int_0^\delta \D_{f(x, u)} V(\varphi(t, x, u)) \spc \mathrm d t
\end{equation}
or, equivalently,
\begin{equation} \label{eqn:optim-based-feedb-optim-2}
	\min_{u \in \U} \abs{ V(\varphi(\delta, x, u)) - V(x)}.
\end{equation}
At every step, $\varphi(\delta, x, u)$, as a solution of \eqref{eqn:sys}, is computed over the sampling time period $[0, \delta]$.
The one step optimization-based control can be defined as
\begin{equation} \label{eqn:optim-based-feedb}
	\kappa_\delta(x) \in \argmin_{u \in \U} V(\varphi(\delta, x, u)).
\end{equation}
This method combines the two steps \eqref{eqn:dini-aiming-direction} and \eqref{eqn:dini-aiming-feeback} in one single optimization problem. 
In comparison to other techniques, $\kappa_\delta(x)$ computed in \eqref{eqn:optim-based-feedb}, explicitly depends on $\delta$.

\subsection{Inf-convolution-based stabilization} \label{sub:infC-stab}

In this technique, a minimizer of \eqref{eqn:InfC} \ie
\begin{equation} \label{eqn:InfC-minimizer}
	y_\alpha(x) \in \arginf_{y \in \R^n} \left\{ V(y) + \frac{1}{2 \alpha^2} \nrm{y - x}^2 \right\}
\end{equation}
is computed to define a proximal subgradient $\zeta_\alpha (x) := \frac{x - y_\alpha(x)}{\alpha^2}$ \citep{Clarke1997-stabilization}.  
A control law is obtained by
\begin{equation} \label{eqn:InfC-ctrl-law}
	\kappa(x) = \arginf_{u \in \U_{\Y}} \scal{\zeta_\alpha(x), f(y_\alpha(x), u)}
\end{equation}
where $\Y$ is a compact set which contains $y_\alpha(x)$.
Numerical studies with the above described methods can be found in \citep{Braun2017-SH-stabilization-Dini-aim,Osinenko2018-practical-SH}.
The next section discusses specifically the use of $F$-disassembled subdifferentials and CLFs for nonsmooth backstepping.

\section{Nonsmooth backstepping} \label{sec:backstepping}

In this section, a variant of nonsmooth backstepping on the example of three-wheel robot with dynamical actuators and dynamical Artstein's circles is presented based on $F$-disassembled subdifferentials and CLFs.

\subsection{Three-wheel robot with dynamical actuators} \label{sub:ENDI}

A three-wheel robot with dynamical control of the driving and steering torques can be described as follows:
\begin{equation} \label{eqn:ENDI}
	\tag {ENDI}
	\begin{aligned}
		\dot x_1 &= \eta_1, & \dot \eta_1 &= u_1,\\
		\dot x_2 &= \eta_2, & \dot \eta_2 &= u_2, \\
		\dot x_3 &= \eta_1 x_2 - x_1 \eta_2. & \\
	\end{aligned}
\end{equation}
The system \eqref{eqn:ENDI}, $\frac{\mathrm d}{\mathrm d t} \begin{bmatrix} x \\ \eta \end{bmatrix} = f_{\text{ENDI}}(x, \eta, u)$, is also called \textit{extended nonholonomic dynamical integrator} (ENDI) \citep{Abbasi2017-backstepping, Sankaranarayanan2009-switched, Pascoal2002-practical}.
It is essentially the Brockett's \textit{nonholonomic integrator} (NI) with additional integrators before the control inputs.
The former reads, accordingly, as 
\begin{equation} \label{eqn:NI}
	\tag {NI}
	\dot x = 
	f_{\text{NI}}(x, u) =	
	\underbrace{\begin{pmatrix}
		1 \\ 0 \\ -x_2
	\end{pmatrix}}_{=: g_1(x)} u_1 + \underbrace{\begin{pmatrix}
		0 \\ 1 \\ x_1
	\end{pmatrix}}_{=: g_2(x)} u_2.
\end{equation}

The following functions are nonsmooth LDGD CLFs for \eqref{eqn:NI} \citep{Braun2017-SH-stabilization-Dini-aim, Clarke2011-discont-stabilization}:
\begin{align*}
	V_1(x) &= x_1^2 + x_2^2 + 2 x_3^2 - 2 \abs{x_3} \sqrt{x_1^2 + x_2^2}, \\
	V_2(x) &= x_1^2 + x_2^2 + 2 x_3^2 + \abs{x_3} (10 - 2 (\abs{x_1} + \abs{x_2})).
\end{align*}
For the ENDI, nonsmooth backstepping based on \citep{Matsumoto2015-position} may be utilized.
To this end, consider the following function \citep{Kimura2015-asymptotic}:
\begin{equation} \label{eqn:NI-F}
	F(x; \theta) := x_1^4 + x_2^4 + \frac{\abs{x_3}^3}{(x_1 \cos (\theta) + x_2 \sin(\theta) + \sqrt{\abs{x_3}})^2}.
\end{equation}

So, $V(x) := \min_{\theta \in \Theta} F(x; \theta)$ is an $F$-disassembled CLF for \eqref{eqn:NI}.
Choose a minimizer $\theta^\star \in \Theta^\star := \argmin_{\theta \in \Theta} F(x; \theta)$ and compute $\zeta(x; \theta) = \nabla_x F(x; \theta)$.
Then, $\zeta(x; \theta^\star) \in \partial_D^F V(x) = \nabla_x F(x; \theta) \mid_{\theta \in \Theta^\star}$ holds.

Use the Sontag's formula \citep{Nakamura2013-asymptotic} to express, formally,

\begin{equation} \label{eqn:feedback-NI}
	\kappa(x; \theta) = -\begin{pmatrix}
		\scal{\zeta(x; \theta), g_1(x)} \\ \scal{\zeta(x; \theta), g_2(x)}
	\end{pmatrix},
\end{equation}
and plugging in a minimizer $\theta^*$ into the above gives a corresponding stabilizing controller for \eqref{eqn:NI}.

We have, for all $x \in \X \sm \{ 0 \}$ and $\theta \in \Theta$,
\begin{equation} \label{eqn:decay-NI-theta}
	\begin{aligned}
		&\scal{\zeta(x; \theta), f_{\text{NI}}(x, u)} \\
		&=\scal{\zeta(x; \theta), g_1(x)} u_1 + \scal{\zeta(x; \theta), g_2(x)} u_2 \\
		&= -\scal{\zeta(x; \theta), g_1(x)}^2 - \scal{\zeta(x; \theta), g_2(x)}^2 < 0.
	\end{aligned}
\end{equation}

Observe, that for all $\theta \in \Theta$, it holds formally that
\begin{equation}
	\scal{\zeta(x; \theta), G(x) \kappa(x; \theta)} \leq 0
\end{equation}
with $G(x) = \begin{bmatrix} g_1(x) & g_2(x)\end{bmatrix} = \begin{pmatrix} 1 & 0 \\ 0 & 1 \\ -x_2 & x_1 \end{pmatrix}$.


Now, the $F$-disassembled CLF \eqref{eqn:NI-F} of \eqref{eqn:NI} is augmented in the spirit of backstepping as follows
\begin{equation} \label{eqn:ENDI-V}
	V_c(x, \eta) = \min_{\theta \in \Theta} \left\{ F(x; \theta) + \frac{1}{2} \nrm{ \eta - \kappa(x; \theta) }^2 \right\}.	
\end{equation}


The $F$-disassembled subdifferential for $\theta_c^\star \in \Theta^\star_c$, where $\Theta^\star_c := \argmin_{\theta \in \Theta} \{ F(x; \theta) + \nicefrac 1 2 \nrm{\eta - \kappa(x; \theta)}^2 \}$ 
is given by
\begin{equation} \label{eqn:disass-grad-ENDI}
	\begin{aligned}
		\partial_D^F V_c(x, \eta) &= \left.\pdiff{F (x; \theta) + \frac{1}{2} \nrm{ \eta - \kappa(x; \theta) }^2 }{\begin{bmatrix} x \\ \eta \end{bmatrix}} \right\vert_{\theta \in \Theta^\star_c} \\
		&= \left. \begin{bmatrix}
		\pdiff{F(x; \theta)}{x} - \pdiff{\kappa(x; \theta)}{x} (\eta - \kappa(x; \theta)) \\
		\eta - \kappa(x; \theta)
		\end{bmatrix} \right \vert_{\theta \in \Theta^\star_c}.
	\end{aligned}
\end{equation}

Define formally (for a generic $\theta$) $\zeta_c(x, \eta; \theta) =$ \\ $ \nabla_{(x,\eta)} F_c(x, \eta; \theta)$ and $z(x;\theta) := \eta - \kappa(x; \theta)$.
Then, express
\begin{equation} \label{eqn:decay-ENDI-theta}
	\begin{aligned}
		&\scal{\zeta_c(x, \eta; \theta), f_{\text{ENDI}}(x, \eta,u)} \\
		&= \scal{		
		\begin{bmatrix}
			\zeta(x; \theta) - \nabla_x \kappa(x; \theta) z \\
			z
		\end{bmatrix},
		\begin{bmatrix}
			G(x) \eta \\ 
			u
		\end{bmatrix}		
		} \\
		&= \scal{\zeta(x; \theta) - \nabla_x \kappa(x; \theta) z, G(x) \eta} + \scal{z, u} \\
		&= \scal{\zeta(x; \theta), G(x) \eta} - \scal{\nabla_x \kappa(x; \theta) z, G(x) \eta} + \scal{z, u} \\
		&= \scal{\zeta(x; \theta), G(x) z} + \scal{\zeta(x; \theta), G(x) \kappa(x; \theta)} \\
		& \ \ \ - \scal{\nabla_x \kappa(x; \theta) z, G(x) \eta} + \scal{z, u} =: S(x, z; \theta). \\
	\end{aligned}
\end{equation}

Let us pick a $\theta^\star_c \in \Theta^\star_c$ and analyze $S(x, z; \theta^\star_c)$.
The following cases are possible.

Case 1: $z(x; \theta^\star_c) = 0$.
Since $z = 0$, $S(x, z; \theta^\star_c)$ reduces to
$\scal{\zeta(x; \theta^\star_c), G(x) \kappa(x; \theta^\star_c)}$.
But, in this case, $\Theta^\star = \Theta^\star_c$ and so the control for \eqref{eqn:NI} is resembled, and $S(x, z; \theta^\star_c) < 0$ accordingly.	

Case 2: $z(x; \theta^\star_c) \not = 0$.
Since
\begin{equation} \label{eqn:feedback-ENDI}
	\begin{aligned}	
		&S(x, z; \theta^\star_c) \\
		&= \scal{z, G(x)^{\mathrm T} \zeta(x; \theta^\star_c)} + \scal{\zeta(x; \theta^\star_c), G(x) \kappa(x; \theta^\star_c)} \\
		& \ \ \ - \scal{z, \nabla_x \kappa(x; \theta^\star_c)^{\mathrm T} G(x) \eta} + \scal{z, u} \\
	\end{aligned}
\end{equation}
holds, a suitable choice for $u$ is \eg
\begin{equation} \label{eqn:ENDI-ctrl}
	\begin{aligned}
		& u = \kappa_c(x, \eta; \theta^\star_c) := \\ 
		& \nabla_x \kappa(x; \theta^\star_c)^{\mathrm T} G(x) \eta - G(x)^{\mathrm T} \zeta(x; \theta^\star_c) - K z
	\end{aligned}
\end{equation}
with $K  > 0$.
This yields
\begin{equation}
	\begin{aligned}
		\scal{\zeta(x; \theta^\star_c), G(x) \kappa(x; \theta^\star_c)} - K \nrm z^2 < 0
	\end{aligned}
\end{equation}
for all $(x,\eta) \in \R^5 \sm \{ 0 \}$.
Thus, the origin of \eqref{eqn:ENDI} is practically stabilized by $\kappa_c(x, \eta; \theta^\star_c)$.
The overall control thus amounts to finding a minimizer $\theta^\star_c$ and plugging it into the expression of $\kappa_c$. 
Notice that $\kappa_c$ contains the expression for $\kappa$ and the corresponding generic form of the subgradient $\zeta$.
As a final note, even a control of the form $u = - K z$ would suffice (in the sense of ultimate upper boundedness of the state) using the fact that the terms in $S$ may be compensated accordingly (details omitted).

\subsection{Dynamical Artstein's circles} \label{sub:Artstein}

The same method as described in the previous section can be applied to the following system, a dynamic extension of the Artstein's Circle:
\begin{equation} \label{eqn:Artsteins-Circle}
	\tag{AC}
	\begin{aligned}
		\dot x_1 &= (-x_1^2 + x_2^2) x_3, \\
		\dot x_2 &= -2 x_1 x_2 x_3, \\
		\dot x_3 &= u, 
	\end{aligned}
\end{equation}
which can be written as
\begin{equation}
	\begin{aligned}
		\dot v & = g(v) w \\
		\dot w &= u
	\end{aligned}
\end{equation}
with $v := \begin{pmatrix} x_1 & x_2 \end{pmatrix}^{\mathrm T} \in \R^2$ and $w := x_3 \in \R$.

A nonsmooth LDGD CLF for $\dot v = g(v) u$ is given as \citep{Braun2018-complete}:
\begin{equation} \label{eqn:CLF-LDGD-Artsteins-Circle}
	V(v) = \sqrt{3 x_1^2 + 4 x_2^2} - \abs{x_1}.
\end{equation}

Referring to \eqref{eqn:NI-F}, the following marginal function gives rise to an $F$-disassembled CLF:
\begin{equation} \label{eqn:CLF-disass-Artsteins-Circle}
	F(v; \theta) = \sqrt 3 x_1 \cos(\theta_1) + 2 x_2 \sin(\theta_1) + x_1 \left( \frac{\theta_2}{\pi} - 1\right).
\end{equation}
Obviously, for \eqref{eqn:CLF-LDGD-Artsteins-Circle} and \eqref{eqn:CLF-disass-Artsteins-Circle}, \eqref{eqn:marginal-fcn} holds with $\Theta = [0, 2 \pi] \times [0, 2 \pi]$.

After computing the $F$-disassembled differential $\partial_D^F V(v) = \nabla_v F(v; \theta) \mid_{\theta \in \Theta^\star}$, choose $\zeta(v; \theta^\star) \in \partial_D^F V(v)$, where $\theta^\star \in \Theta^\star := \argmin_{\theta \in \Theta} F(v; \theta)$.
Then, with $\kappa(v; \theta) = - \scal{ \zeta(v; \theta), g(v) }$ as a feedback, it follows
\begin{equation}
	\scal{ \zeta(v; \theta), g(v) \kappa(v; \theta)} = -\scal{ \zeta(v; \theta), g(v) }^2 < 0
\end{equation}
which holds also for all $\theta \in \Theta$ and all $v \in \R^2 \sm \{ 0 \}$, and practically stabilizes $\dot v = g(v)w$.
Now, \eqref{eqn:CLF-disass-Artsteins-Circle} can be extended via backstepping to
\begin{equation}
	V_c(v, w) = \min_{\theta \in \Theta_c} \left\{ F(v; \theta) + \frac 1 2 \nrm{ w - \kappa(v; \theta) }^2 \right\}.	
\end{equation}
Define $z := w - \kappa(v; \theta)$.
Since the disassembled subdifferential is given similarly to \eqref{eqn:disass-grad-ENDI}, one can show, that for $\zeta_c(v, w; \theta^\star) \in \partial_D^F V_C(v, w)$ with $\theta^\star \in \Theta^\star_c$, $\Theta^\star_c$ defined as in Section \ref{sub:ENDI}, the following equality holds:
\begin{equation}
	\begin{aligned}
		&\scal{		
		\begin{bmatrix}
			\zeta(v; \theta) - \nabla_v \kappa(v; \theta) z \\
			z
		\end{bmatrix},
		\begin{bmatrix}
			g w \\ 
			u
		\end{bmatrix}		
		} \\
		&= \scal{\zeta(v; \theta) - \nabla_v \kappa(v; \theta) z, g w} + \scal{z, u} \\
		&= \scal{\zeta(v; \theta), g w} - \scal{\nabla_v \kappa(v; \theta) z, g w} + \scal{z, u} \\
		&= \scal{\zeta(v; \theta), g z} + \scal{\zeta(v; \theta), g \kappa(v; \theta)} \\
		& \ \ \ - \scal{\nabla_v \kappa(v; \theta) z, g w} + \scal{z, u}. \\
	\end{aligned}
\end{equation}
The corresponding control law can be derived in a similar manner as in Section \ref{sub:ENDI}.

%

\section{Sliding mode} \label{sec:SMC}



Consider a system
\begin{equation}
\label{eqn:sys-SMC}
\begin{array}{ll}
\dot x & = f(x, u), x \in \R^n, \\
u & = \kappa(x) + s(x) \sigma(x),
\end{array}
\end{equation}
where $\kappa, s$ are continuous and
\begin{equation}
\label{eqn:switch-ctrl}
\sigma(x) = \begin{cases}
1, & \chi(x) > 0, \\
-1, & \chi(x) < 0
\end{cases}
\end{equation}
describes the discontinuous part of the controller around a sliding surface $\Sigma = \{x: \chi(x) = 0 \}$.
To process \eqref{eqn:sys-SMC}, unlike in the usual case, two Lyapunov-like functions $V_1, V_2$ are employed that satisfy, for some $\bar w_1 > 0$,
\begin{align*}
& \hspace{-30pt} \begin{array}{ll}
& V_1 \text{ and } V_1 + V_2 \text{ are proper}, V_1(x) = 0 \iff x \in \Sigma, \\  
& \text{$V_1$ is continuously differentiable on } \R^n \setminus \Sigma \text{ and} \\
& \langle \nabla V_1(x) , f(x, \kappa(x) + s(x)\sigma(x))\rangle \le - \bar w_1, x \in \R^n \setminus \Sigma;
\end{array} \\
& \hspace{-30pt} \begin{array}{ll}
& V_2(0) = 0, w_2(0) = 0, \\
& x \in \Sigma \setminus \{0\} \implies V_2(x) > 0, w_2(x) > 0, \\
& \text{$V_2$ is continuously differentiable on } \R^n \setminus \{0\} \text{ and} \\
& \sup_{v \in \partial_L^f(x) } \langle \nabla V_2(x), v \rangle \le -w_2(x), x \in \Sigma \setminus \{0\},
\end{array}
\end{align*}
where $\partial_L^f(x) := \{ \lim_{i \ra \infty} f(x, \kappa(x) + s(x) \sigma_i) : x_i \ra x, \sigma_i \in \sigma(x_i) \}$ describes the limiting behavior of possible system velocities (cf. the construction in Definition \ref{dfn:limit-subdiff}). 
Let $f_1(x) := f(x, \kappa(x) + s(x))$ and $f_2(x) := f(x, \kappa(x) - s(x))$, and $\X_1 := \{ \chi(x) > 0 \}, \X_2 := \{ \chi(x) < 0 \}$.
In general, \eqref{eqn:sys-SMC} in SMC mode is usually treated in the sense of a differential inclusion
\citep{Perruquetti2002-slid-mode,Slotine1991-nonlin-ctrl}:
\begin{equation}
\label{eqn:SMC-DI}
\dot x	\in F(x) = \begin{cases}
f_1(x), & x	\in \X_1,\\
(1 - \alpha)f1(x) + \alpha f_2(x), \alpha \in [0, 1], & x \in \Sigma, \\
f_2(x), & x	\in \X_2.
\end{cases}
\end{equation} 
This is a particular example of \textit{Filippov regularization} \citep{Cortes2008-discont-dyn-sys} and, in this case, is upper semi-continuous, has closed and convex images and thus admits a Filippov solution \citep{Zabczyk2009-mathematical}, which is an absolutely continuous function whose derivative satisfies the differential inclusion almost everywhere.
A Filippov solution is an idealized construction describing a perfect sliding mode in this context.
In practice, since the control is usually realized digitally, S\&H analysis of SMC comes in handy \citep{Clarke2009-slid-mode-stab}.
The general idea thereof is to chose a sampling time bound small enough that the system approach a specially chosen vicinity of $\Sigma$ where $V_2$ still retains some of it decay rate $w_2$ (by a continuity argument).
The attraction to this vicinity is in turn ensured by retaining some of the decay rate $\bar w_1$.
Combining the two decay properties together ensures practical stabilization of \eqref{eqn:sys-SMC} in the S\&H mode \citep{Clarke2009-slid-mode-stab}.

\section{Computational aspects} \label{sec:comp-aspects}

The described practical stabilization methods of the previous sections rely on computation of primarily two things: \begin{enumerate*}
	\item allowed sampling time for the desired stabilization precision;
	\item control actions.
\end{enumerate*}
General techniques of Section \ref{sec:stabilization} rely heavily on various optimizations.
It was shown that optimization accuracy greatly influences stabilization precision \citep{Osinenko2018-practical-SH}.
Moreover, the involved LDGD CLF must satisfy certain regularity properties.
In particular, the following:
For all compact sets $\Y, \F \subset \R^n$ and for all $\nu, \chi > 0$ there exist $\tilde \Y \subseteq \Y, \mu \ge 0$ such that:
\begin{enumerate}
	\item For each $\tilde y \in  \tilde \Y, \theta \in \F$ and $\forall \mu' \in (0, \mu]$ it holds that
	\begin{equation} \label{eqn-loc-hom-cond}
	\tag{hom} 
	\abs{ \frac{ V(\tilde y + \mu' \theta) - V(\tilde y) }{\mu'} - \D_\theta V(\tilde y) } \le \nu;
	\end{equation}
	\item For each $y \in \Y$ there exists $\tilde y \in \tilde  \Y$ such that
	\begin{equation} \label{eqn-neighb-pt-cond}
	\tag{npt} 
	\nrm{y - \tilde y} \le \chi.
	\end{equation}
\end{enumerate} 
Under these conditions, one can find bounds on the optimization accuracy of the inf-convolution-based method so as to achieve the desired stabilization precision.
Furthermore, it was shown in \citep{Osinenko2018-constr-SMC} what regularity properties of the involved Lyapunov-like functions have to be satisfied to effectively compute required bounds on the sampling time.
The respective machinery was addressed on the example of SMC.
In contrast to \citep{Clarke2009-slid-mode-stab}, an actual numerical example of practical SMC stabilization for vehicle slip control was shown.
The computed sampling time bounds were satisfactory for the considered application.

Here, we give a short case study that demonstrates the effect of the described computational uncertainty.
Namely, for \eqref{eqn:ENDI} and its CLF \eqref{eqn:ENDI-V}, a control law was computed in the S\&H framework using \eqref{eqn:InfC-minimizer} and \eqref{eqn:InfC-ctrl-law} with optimization accuracy $\eps$ and $\gamma$, respectively.
The initial condition is set to $x_0 = \begin{pmatrix} -1 & 0.5 & 0.01 & 0.05 & 0.075 \end{pmatrix}^{\mathrm T}$ and the set of admissible controls is given as $\U = [-3,3]$.
Furthermore, we set $\alpha$ in \eqref{eqn:InfC} to $\alpha = 0.1$ and the sampling time $\delta$ to $\delta = 0.005$. \linebreak

\begin{figure}
\centering
  \begin{tikzpicture}
	\coordinate (1) at (0,0);
  	\node [right=of 1, scale = 0.8] {\includegraphics[width = 0.5\textwidth]{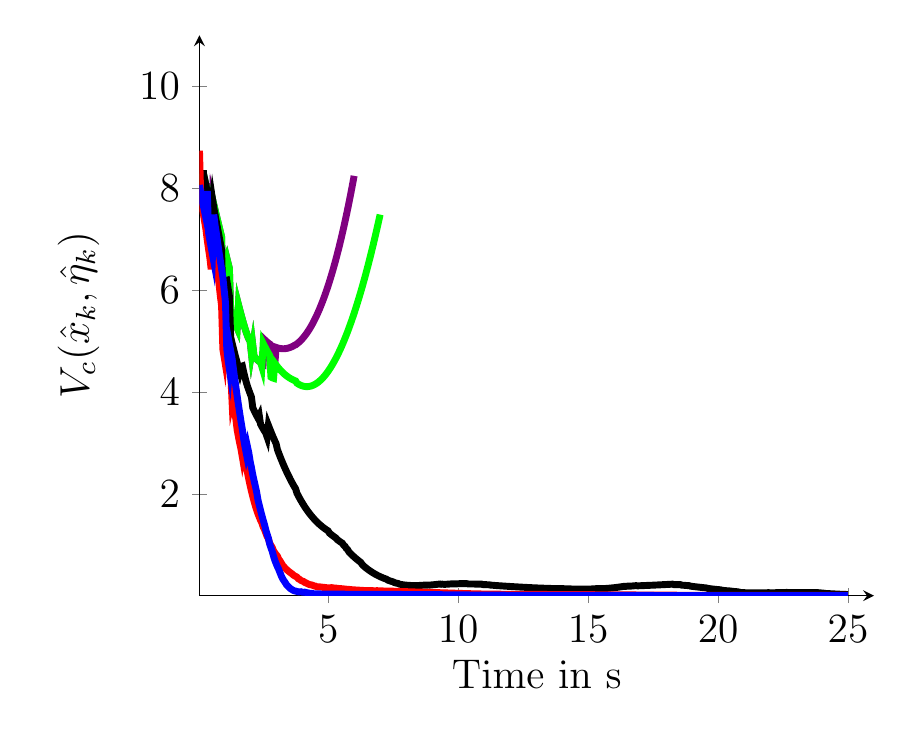}};	
  	\coordinate (2) at (1.875,1.75);
  	\node [right=of 2, scale = 0.8] {\includegraphics[width = 0.5\textwidth]{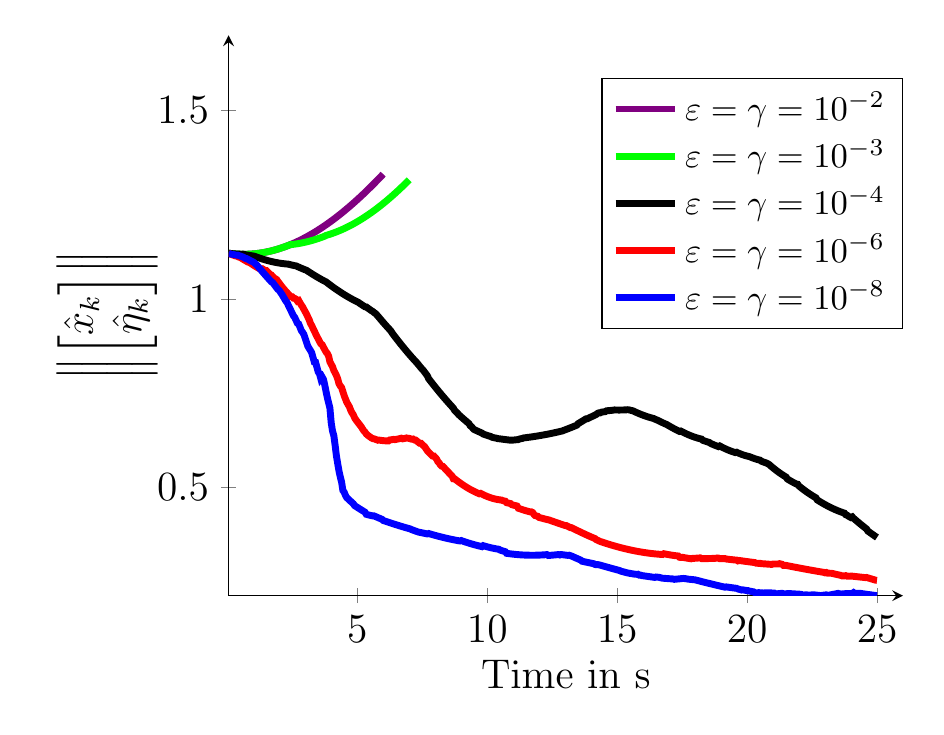}};
	\end{tikzpicture}
	\caption{Practical stabilization of the three-wheel robot with dynamic actuators under different computational accuracy.}
	\label{fig:gfx_epsilon}
\end{figure}

The influence of the optimization accuracy on the state and CLF behavior can be seen in Fig. \ref{fig:gfx_epsilon} for different values of $\eps$ and $\gamma$, namely $\eps = \gamma \in \{ 10^{-2}, 10^{-3}, 10^{-4}, 10^{-6}, 10^{-8} \}$.
It can be observed that insufficient accuracy ($\eps = \gamma = 10^{-3}$ and $\eps = \gamma = 10^{-2}$) leads instability.
Higher accuracies lead to ever smaller vicinities of the origin that the state converges into.
This clearly demonstrates that computational uncertainty must be taken into account in practical stabilization.

\section{Conclusion}

This work surveyed, in a brief form, some of the key modern nonsmooth stabilization tools and techniques.
These include general practical stabilization, sliding-mode control, nonsmooth backstepping.
Examples were provided.
In addition, this work briefly discussed some computational aspects of practical stabilization one should be concerned about in applications.

\bibliography{bib/constructing-LFs,bib/discont-DE,bib/dyn-sys,bib/MPC,bib/non-smooth-analysis,bib/opt-ctrl,bib/sliding-mode,bib/stability,bib/stabilization,bib/nonlin-ctrl,bib/semiconcave,bib/PDE,bib/Osinenko,bib/Filippov-sol}                                                                                  

\begin{thebibliography}{42}
\providecommand{\natexlab}[1]{#1}
\providecommand{\url}[1]{\texttt{#1}}
\providecommand{\urlprefix}{URL }
\expandafter\ifx\csname urlstyle\endcsname\relax
  \providecommand{\doi}[1]{doi:\discretionary{}{}{}#1}\else
  \providecommand{\doi}{doi:\discretionary{}{}{}\begingroup
  \urlstyle{rm}\Url}\fi

\bibitem[{Abbasi et~al.(2017)Abbasi, ur~Rehman, and
  Shah}]{Abbasi2017-backstepping}
Abbasi, W., ur~Rehman, F., and Shah, I. (2017).
\newblock Backstepping based nonlinear adaptive control for the extended
  nonholonomic double integrator.
\newblock \emph{Kybernetika}, 53(4), 578--594.

\bibitem[{Aeyels(1985)}]{Aeyels1985-stabilization}
Aeyels, D. (1985).
\newblock Stabilization of a class of nonlinear systems by a smooth feedback
  control.
\newblock \emph{Systems \& Control Letters}, 5(5), 289--294.

\bibitem[{Braun et~al.(2017)Braun, Gr{\"u}ne, and
  Kellett}]{Braun2017-SH-stabilization-Dini-aim}
Braun, P., Gr{\"u}ne, L., and Kellett, C. (2017).
\newblock Feedback design using nonsmooth control {L}yapunov functions: A
  numerical case study for the nonholonomic integrator.
\newblock In \emph{Proceedings of the 56th IEEE Conference on Decision and
  Control}.

\bibitem[{Braun et~al.(2018)Braun, Gr{\"u}ne, and Kellett}]{Braun2018-complete}
Braun, P., Gr{\"u}ne, L., and Kellett, C. (2018).
\newblock Complete instability of differential inclusions using {L}yapunov
  methods.
\newblock In \emph{2018 IEEE Conference on Decision and Control (CDC)},
  718--724. IEEE.

\bibitem[{Brockett(1983)}]{Brockett1983-stabilization}
Brockett, R. (1983).
\newblock Asymptotic stability and feedback stabilization.
\newblock \emph{Differential geometric control theory}, 27(1), 181--191.

\bibitem[{Camilli et~al.(2008)Camilli, Gr{\"u}ne, and
  Wirth}]{Camilli2008-CLF-by-Zubov-mtd}
Camilli, F., Gr{\"u}ne, L., and Wirth, F. (2008).
\newblock Control {L}yapunov functions and {Z}ubov's method.
\newblock \emph{SIAM Journal on Control and Optimization}, 47(1), 301--326.

\bibitem[{Cannarsa and Sinestrari(2004)}]{Cannarsa2004-semiconcave}
Cannarsa, P. and Sinestrari, C. (2004).
\newblock \emph{Semiconcave functions, Hamilton-Jacobi equations, and optimal
  control}, volume~58.
\newblock Springer Science \& Business Media.

\bibitem[{Clarke(2011)}]{Clarke2011-discont-stabilization}
Clarke, F. (2011).
\newblock Lyapunov functions and discontinuous stabilizing feedback.
\newblock \emph{{A}nnual {R}eviews in {C}ontrol}, 35(1), 13--33.

\bibitem[{Clarke et~al.(1997)Clarke, Ledyaev, Sontag, and
  Subbotin}]{Clarke1997-stabilization}
Clarke, F., Ledyaev, Y., Sontag, E., and Subbotin, A. (1997).
\newblock Asymptotic controllability implies feedback stabilization.
\newblock \emph{IEEE Transactions on Automatic Control}, 42(10), 1394--1407.

\bibitem[{Clarke et~al.(1995)Clarke, Ledyaev, Stern, and
  Wolenski}]{Clarke1995-sys-traj}
Clarke, F., Ledyaev, Y., Stern, R., and Wolenski, P. (1995).
\newblock Qualitative properties of trajectories of control systems: a survey.
\newblock \emph{Journal of dynamical and control systems}, 1(1), 1--48.

\bibitem[{Clarke et~al.(2008)Clarke, Ledyaev, Stern, and
  Wolenski}]{Clarke2008-nonsmooth-analys}
Clarke, F., Ledyaev, Y., Stern, R., and Wolenski, P. (2008).
\newblock \emph{Nonsmooth Analysis and Control Theory}, volume 178.
\newblock Springer Science \& Business Media.

\bibitem[{Clarke and Vinter(2009)}]{Clarke2009-slid-mode-stab}
Clarke, F. and Vinter, R. (2009).
\newblock Stability analysis of sliding-mode feedback control.
\newblock \emph{Control and Cybernetics}, 4(38), 1169--1192.

\bibitem[{Coron(1992)}]{Coron1992-stabilization}
Coron, J.M. (1992).
\newblock Global asymptotic stabilization for controllable systems without
  drift.
\newblock \emph{Mathematics of Control, Signals and Systems}, 5(3), 295--312.

\bibitem[{Coron(1995)}]{Coron1995-stabilization}
Coron, J.M. (1995).
\newblock On the stabilization in finite time of locally controllable systems
  by means of continuous time-varying feedback law.
\newblock \emph{SIAM Journal on Control and Optimization}, 33(3), 804--833.

\bibitem[{Coron and d'Andrea Novel(1991)}]{Coron1991-smooth-stabilization}
Coron, J.M. and d'Andrea Novel, B. (1991).
\newblock Smooth stabilizing time-varying control laws for a class of nonlinear
  systems. {A}pplication to mobile robots.
\newblock In \emph{Nonlinear Control Systems Design 1992. Selected Papers from
  the 2nd IFAC Symposium}, 413--18.

\bibitem[{Coron and Pomet(1993)}]{Coron1993-stab-time-var}
Coron, J.M. and Pomet, J.B. (1993).
\newblock A remark on the design of time-varying stabilizing feedback laws for
  controllable systems without drift.
\newblock In \emph{Nonlinear Control Systems Design 1992}, 397--401. Elsevier.

\bibitem[{Coron and Rosier(1994)}]{Coron1994-stabilization}
Coron, J.M. and Rosier, L. (1994).
\newblock A relation between continuous time-varying and discontinuous feedback
  stabilization.
\newblock \emph{Journal of Mathematical Systems, Estimation, and Control}, 4,
  67--84.

\bibitem[{Cortes(2008)}]{Cortes2008-discont-dyn-sys}
Cortes, J. (2008).
\newblock Discontinuous dynamical systems.
\newblock \emph{IEEE Control Systems Magazine}, 28(3), 36--73.

\bibitem[{Crandall and Lions(1983)}]{Crandall1983-viscosity-HJB}
Crandall, M. and Lions, P.L. (1983).
\newblock Viscosity solutions of {H}amilton-{J}acobi equations.
\newblock \emph{Transactions of the American mathematical society}, 277(1),
  1--42.

\bibitem[{Filippov(1988)}]{Filippov2013-discont-dyn-sys}
Filippov, A. (1988).
\newblock \emph{Differential {E}quations with {D}iscontinuous {R}ighthand
  {S}ides: {C}ontrol {S}ystems}.
\newblock Springer Science \& Business Media.

\bibitem[{Fridland and Levant(1999)}]{Fridland1999-slid-mode}
Fridland, L. and Levant, A. (1999).
\newblock \emph{{H}igher {O}rder {S}liding {M}odes}.
\newblock Sliding mode in Automatic Control, Ecole Central de Lille.

\bibitem[{Fridman and Levant(1996)}]{Fridman1996-slid-mode}
Fridman, L. and Levant, A. (1996).
\newblock Higher order sliding modes as a natural phenomenon in control theory.
\newblock In \emph{Robust Control via variable structure and Lyapunov
  techniques}, 107--133. Springer.

\bibitem[{Kawski(1989)}]{Kawski1989-stabilization-plane}
Kawski, M. (1989).
\newblock Stabilization of nonlinear systems in the plane.
\newblock \emph{Systems \& Control Letters}, 12(2), 169--175.

\bibitem[{Kellett et~al.(2004)Kellett, Shim, and Teel}]{Kellett2004-Dini-aim}
Kellett, C., Shim, H., and Teel, A. (2004).
\newblock Further results on robustness of (possibly discontinuous) sample and
  hold feedback.
\newblock \emph{IEEE Transactions on Automatic Control}, 49(7), 1081--1089.

\bibitem[{Kellett and Teel(2000)}]{Kellett2000-Dini-aim}
Kellett, C. and Teel, A. (2000).
\newblock Uniform asymptotic controllability to a set implies locally
  {L}ipschitz control-{L}yapunov function.
\newblock In \emph{Proceedings of the 39th IEEE Conference on Decision and
  Control}, volume~4, 3994--3999.

\bibitem[{Khaneja and Brockett(1999)}]{Khaneja1999-dynamic-feedback}
Khaneja, N. and Brockett, R. (1999).
\newblock Dynamic feedback stabilization of nonholonomic systems (i).
\newblock In \emph{IEEE Conference on Decision and Control}, volume~2,
  1640--1645. IEEE.

\bibitem[{Kimura et~al.(2015)Kimura, Nakamura, and
  Yamashita}]{Kimura2015-asymptotic}
Kimura, S., Nakamura, H., and Yamashita, Y. (2015).
\newblock Asymptotic stabilization of two-wheeled mobile robot via locally
  semiconcave generalized homogeneous control {L}yapunov function.
\newblock \emph{SICE Journal of Control, Measurement, and System Integration},
  8(2), 122--130.

\bibitem[{Matsumoto et~al.(2015)Matsumoto, Nakamura, Satoh, and
  Kimura}]{Matsumoto2015-position}
Matsumoto, R., Nakamura, H., Satoh, Y., and Kimura, S. (2015).
\newblock Position control of two-wheeled mobile robot via semiconcave function
  backstepping.
\newblock In \emph{2015 IEEE Conference on Control Applications (CCA)},
  882--887. IEEE.

\bibitem[{Morin et~al.(1999)Morin, Pomet, and Samson}]{Morin1999-design}
Morin, P., Pomet, J., and Samson, C. (1999).
\newblock Design of homogeneous time-varying stabilizing control laws for
  driftless controllable systems via oscillatory approximation of lie brackets
  in closed loop.
\newblock \emph{SIAM Journal on Control and Optimization}, 38(1), 22--49.

\bibitem[{Nakamura et~al.(2013)Nakamura, Tsuzuki, Fukui, and
  Nakamura}]{Nakamura2013-asymptotic}
Nakamura, H., Tsuzuki, T., Fukui, Y., and Nakamura, N. (2013).
\newblock Asymptotic stabilization with locally semiconcave control {L}yapunov
  functions on general manifolds.
\newblock \emph{Systems \& Control Letters}, 62(10), 902--909.

\bibitem[{Osinenko et~al.(2018{\natexlab{a}})Osinenko, Beckenbach, and
  Streif}]{Osinenko2018-practical-SH}
Osinenko, P., Beckenbach, L., and Streif, S. (2018{\natexlab{a}}).
\newblock Practical sample-and-hold stabilization of nonlinear systems under
  approximate optimizers.
\newblock \emph{IEEE Control Systems Letters (L-CSS)}, 2(4), 569--574.

\bibitem[{Osinenko et~al.(2018{\natexlab{b}})Osinenko, Devadze, and
  Streif}]{Osinenko2018-constr-SMC}
Osinenko, P., Devadze, G., and Streif, S. (2018{\natexlab{b}}).
\newblock Practical stability analysis of sliding-mode control with explicit
  computation of sampling time.
\newblock \emph{Asian Journal of Control}.

\bibitem[{Pascoal and Aguiar(2002)}]{Pascoal2002-practical}
Pascoal, A. and Aguiar, A. (2002).
\newblock Practical stabilization of the extended nonholonomic double
  integrator.
\newblock \emph{Proc. 10th Mediterranean Conferenceon Control and Automation}.

\bibitem[{Perruquetti and Barbot(2002)}]{Perruquetti2002-slid-mode}
Perruquetti, W. and Barbot, J.P. (2002).
\newblock \emph{{S}liding {M}ode {C}ontrol in {E}ngineering}.
\newblock CRC Press.

\bibitem[{Pomet(1992)}]{Pomet1992-dyn-stabilization}
Pomet, J.B. (1992).
\newblock Explicit design of time-varying stabilizing control laws for a class
  of controllable systems without drift.
\newblock \emph{Systems \& Control Letters}, 18(2), 147--158.

\bibitem[{Ryan(1994)}]{Ryan1994-stabilization}
Ryan, E. (1994).
\newblock On {B}rockett's condition for smooth stabilizability and its
  necessity in a context of nonsmooth feedback.
\newblock \emph{SIAM Journal on Control and Optimization}, 32(6), 1597--1604.

\bibitem[{Samson(1991)}]{Samson1991-stabilization}
Samson, C. (1991).
\newblock Velocity and torque feedback control of a nonholonomic cart.
\newblock In \emph{Advanced robot control}, 125--151. Springer.

\bibitem[{Sankaranarayanan and
  Mahindrakar(2009)}]{Sankaranarayanan2009-switched}
Sankaranarayanan, V. and Mahindrakar, A. (2009).
\newblock Switched control of a nonholonomic mobile robot.
\newblock \emph{Communications in Nonlinear Science and Numerical Simulation},
  14(5), 2319--2327.

\bibitem[{Slotine and Li(1991)}]{Slotine1991-nonlin-ctrl}
Slotine, J.J.E. and Li, W.A. (1991).
\newblock \emph{{A}pplied {N}onlinear {C}ontrol}.
\newblock Prentice Hall.

\bibitem[{Subbotin(2013)}]{Subbotin2013-gen-sol-PDE}
Subbotin, A.I. (2013).
\newblock \emph{Generalized Solutions of First-Order PDEs: the Dynamical
  Optimization Perspective}.
\newblock Springer Science \& Business Media.

\bibitem[{Young et~al.(1996)Young, Utkin, and Ozguner}]{Young1996-slid-mode}
Young, K., Utkin, V.I., and Ozguner, U. (1996).
\newblock A control engineer's guide to sliding mode control.
\newblock In \emph{Proceedings of the IEEE International Workshop on Variable
  Structure Systems VSS'96}, 1--14.

\bibitem[{Zabczyk(2009)}]{Zabczyk2009-mathematical}
Zabczyk, J. (2009).
\newblock \emph{Mathematical control theory: an introduction}.
\newblock Springer Science \& Business Media.

\end{thebibliography}

\end{document}